\documentclass{article}
\usepackage{authblk}
\usepackage{amsmath,amsfonts}
\usepackage{array}
\usepackage[caption=false,font=normalsize,labelfont=sf,textfont=sf]{subfig}
\usepackage{textcomp}
\usepackage{stfloats}
\usepackage{hyperref}
\usepackage{url}
\usepackage{soul}
\usepackage{verbatim}
\usepackage{cite}
\usepackage{xcolor}
\usepackage{graphicx}
\usepackage{tikz}
\usetikzlibrary{fit,calc}
\def\BibTeX{{\rm B\kern-.05em{\sc i\kern-.025em b}\kern-.08em
    T\kern-.1667em\lower.7ex\hbox{E}\kern-.125emX}}
\usepackage{balance}
\usepackage{url}

\makeatletter
\newcommand\email[2][]%
   {\newaffiltrue\let\AB@blk@and\AB@pand
      \if\relax#1\relax\def\AB@note{\AB@thenote}\else\def\AB@note{\relax}%
        \setcounter{Maxaffil}{0}\fi
      \begingroup
        \let\protect\@unexpandable@protect
        \def\thanks{\protect\thanks}\def\footnote{\protect\footnote}%
        \@temptokena=\expandafter{\AB@authors}%
        {\def\\{\protect\\\protect\Affilfont}\xdef\AB@temp{#2}}%
         \xdef\AB@authors{\the\@temptokena\AB@las\AB@au@str
         \protect\\[\affilsep]\protect\Affilfont\AB@temp}%
         \gdef\AB@las{}\gdef\AB@au@str{}%
        {\def\\{, \ignorespaces}\xdef\AB@temp{#2}}%
        \@temptokena=\expandafter{\AB@affillist}%
        \xdef\AB@affillist{\the\@temptokena \AB@affilsep
          \AB@affilnote{}\protect\Affilfont\AB@temp}%
      \endgroup
       \let\AB@affilsep\AB@affilsepx
}
\makeatother

\renewcommand\Affilfont{\fontsize{9}{10.8}\itshape}

\title{Reduced Order Modeling for Real-Time Monitoring of Structural Displacements due to Electromagnetic Forces in Large Scale Tokamaks}
\author[1]{Francesco Lucchini}
\author[3]{Alessandro Frescura}
\author[1,2]{Riccardo Torchio}
\author[1]{Piergiorgio Alotto}
\author[1,3]{Paolo Bettini}
\affil[1]{Department of Industrial Engineering, University of Padova, Via Gradenigo 6/a, 35131, Padova, Italy}
\email{francesco.lucchini@unipd.it}
\affil[2]{Department of Information Engineering, University of Padova, Via Gradenigo 6/b, 35131, Padova, Italy}
\affil[3]{CRF, University of Padova, Padova, Italy}
\date{} 

\begin{document}

\maketitle






\begin{abstract}
The real-time monitoring of the structural displacement of the Vacuum Vessel (VV) of thermonuclear fusion devices caused by electromagnetic (EM) loads is of great interest. 
In this paper, Model Order Reduction (MOR) is applied to the Integral Equation Methods (IEM) and the Finite Elements  Method (FEM)
to develop Electromagnetic and Structural Reduced Order Models (ROMs) compatible with real-time execution which allows for the real-time monitoring of strain and displacement in critical positions of Tokamaks machines.
Low-rank compression techniques based on hierarchical matrices are applied to reduce the computational cost during the offline stage when the ROMs are constructed.  Numerical results show the accuracy of the approach and demonstrate the compatibility with real-time execution in standard hardware.  
\end{abstract}


\section{Introduction}\label{sec.intro}
The evaluation of the structural deformation of conductive structures in magnetic confinement fusion devices is of great interest for the operation of current experimental and future machines. Structural deformation is generally due to loads of different nature, for example thermal, seismic, and electromagnetic (EM)~\cite{takeda2004design}.

In pulsed-operated Tokamaks such as ITER (International Thermonuclear Experimental Reactor), abrupt plasma terminations may result in fast disruption events in the form of Major Disruptions (MD)~\cite{boozer2012theory} or symmetric and asymmetric Vertical Displacement Events (VDE). The latter, occurring in elongated plasmas, can be traced back to the failure of the plasma vertical position control system~\cite{paccagnella2005vertical,hassanein2008vertical}. During these occurrences, the eddy currents induced in passive conductive structures, such as the Vacuum Vessel (VV), are accountable for generating electromagnetic stress through the Lorentz force, resulting in structural displacements. 
It is important to highlight that modeling EM forces serves not only to estimate the EM stresses resulting from the mentioned phenomena but also holds significant relevance during regular operational conditions~\cite{albanese2010electromagnetic,albanese2015effects,pustovitov2017computation,yanovskiy2022sideways,sadakov2024practical}. 

In realistic applications, the monitoring of the mechanical stress on conductive structures caused by EM forces typically relies on a series of sensors, i.e.,  Strain Gauges (SGs)~\cite{6635411}. However, the number of SGs that can be installed in the machine is limited and they cannot be placed in all the points of interest. Thus, the possibility of exploiting a Digital Twin, i.e., a virtual model of the machine serving as a virtual sensor of the displacements, is very attractive.  Naturally, the computational cost of the model must remain low to enable real-time execution. It is worth noting that the notion of Digital Twins has newly surfaced within the realm of thermonuclear fusion~\cite{iglesias2017digital,kwon2022development}.
Obviously, to construct a virtual EM model of the machine it is mandatory to also have a real-time model of the plasma current. Concerning this issue, the interested reader is referred to the
existing literature about this topic~\cite{felici2011real}.

In this paper, a combined electromagnetic-structural  Model Order Reduction (MOR)  for the real-time monitoring of mechanical stress on the conductive parts of thermonuclear fusion devices due to the electromagnetic loads caused by VDE, is proposed. The EM problem can be tackled by using both the Finite Element Method (FEM) and the Volume Integral Equation (VIE) method. In this paper, a VIE method is adopted following the existing literature \cite{albanese1988integral,voltolina2021optimized}. However, the proposed methodology can be similarly applied to FEM. 
To solve the problem of dense matrices generated by the VIE, low-rank compression techniques based on hierarchical matrices ($\mathcal{H}$-matrices) are used. 
The structural problem is instead solved by using the FEM. 
The Proper Orthogonal Decomposition (POD) method is employed to develop Reduced Order Models (ROMs) for both the EM and structural problems, with special attention devoted to minimizing the computational expense associated with coupling the two physics.

The remainder of the article is organized as follows. In section~\ref{sec.VIE}, the VIE method for the analysis of the eddy currents is introduced and a data-sparse representation based on hierarchical matrices is applied to efficiently handle the fully populated matrices. Section~\ref{sec.FEM} describes the linear elasticity model for the structural analysis, while the coupling with the EM problem is delineated in section~\ref{sec.coupling_VIE_FEM}. The  MOR for the EM and structural problems is detailed in section~\ref{sec.pMOR}. A brief description of the steps for the solution of the time-dependent problem is given in section~\ref{sec.time_integral}, while the numerical results are outlined in section~\ref{sec.numeric_ex} and the conclusions are finally given in section~\ref{sec.conclusion}.

\section{Volume integral method for solving eddy currents}\label{sec.VIE}
A common starting point for the solution of the time-domain eddy currents problems in $[0,T]\times\Omega$ where $\Omega\subset\mathbb{R}^3$ and $[0,T]$ is the temporal range, is the Electric Field Integral Equation (EFIE) coupled with the continuity relation~\cite{passarotto2022foundations}:
\begin{align}
&\mathbf{E}(t,\mathbf{r})=-\frac{\partial \mathbf{A}(t,\mathbf{r})}{\partial t}-\nabla\varphi(t,\mathbf{r}) - \frac{\partial\mathbf{A}_{ext}(t,\mathbf{r})}{\partial t}
\nonumber \\
&\nabla\cdot \mathbf{J}(t,\mathbf{r}) = 0 \label{eq.efie_div},
\end{align}
where $\mathbf{E}(t,\mathbf{r})$ is the electric field, $\mathbf{J}(t,\mathbf{r})$ is the current density, $\varphi(t,\mathbf{r})$ the electric scalar potential, and $\mathbf{A}_{ext}(t,\mathbf{r})$ is the magnetic vector potential related, for example, to external excitations. Note that inside the conductive domain,  $\mathbf{E}(t,\mathbf{r})$ and $\mathbf{J}(t,\mathbf{r})$ are connected through the constitutive relation:
\begin{equation}\label{eq.EJ_const}
    \mathbf{E}(t,\mathbf{r})=\rho\mathbf{J}(t,\mathbf{r}),
\end{equation}
where $\rho$ is the electric resistivity.
The magnetic vector potential $\mathbf{A}(t,\mathbf{r})$ can be written in terms of the current density field by introducing the static Green's function: 
\begin{equation} \label{eq.Green}
G(\mathbf{r},\mathbf{r}^\prime)=\frac{1}{4\pi\lVert\mathbf{r}-\mathbf{r}^\prime\rVert},
\end{equation}
as:
\begin{equation} \label{eq.A}
\mathbf{A}(t,\mathbf{r})=\mu_0\int_{\Omega}G(\mathbf{r},\mathbf{r}^\prime)\mathbf{J}(t,\mathbf{r}^\prime)d\Omega,
\end{equation} 
where $\mathbf{r}$ identifies the target point and $\mathbf{r}^\prime$ is the source point.
Historically, challenges posed by eddy currents were addressed through the use of the electric vector potential formulation. This approach leverages the divergence-free condition of the current density through the definition of the electric vector potential $\mathbf{T}$ such that $\mathbf{J}=\nabla \times \mathbf{T}$~\cite{albanese1988integral,albanese2008coupling,bettini2013numerical,bettini2014computation}. The $\mathbf{T}$-formulation has the advantage of reducing the number of unknowns of the problem, however, if the computational domain $\Omega$ is multiply connected, as is usually the case for the vacuum vessels of magnetic confinement fusion devices, the concept of cohomology must be introduced~\cite{voltolina2021optimized,bettini2014lazy}. To avoid the generation of the cohomology group, here the magneto quasi-static (MQS) EFIE \eqref{eq.efie_div}, has been solved with a $\mathbf{J}-\varphi$ formulation, where both the current density and the scalar potential are unknowns. Starting from the discretization of the computational domain $\Omega$ with $N_v$ elements for a total of $N_f$ faces, the current density and electric potential are expanded in terms of vector, $\mathbf{w}(\mathbf{r})$, and piecewise constant, $\psi(\mathbf{r})$, basis functions as: 
\begin{equation}\label{eq.j_phi_basis}
    \begin{array}{cc}
        \mathbf{J}(t,\mathbf{r})=\sum_{k=1}^{N_f}j(t)\mathbf{w}(\mathbf{r}), & 
        \Phi(t,\mathbf{r})=\sum_{h=1}^{N_v}\phi(t)\psi(\mathbf{r}),
    \end{array}
\end{equation}
where $j(t)$ and $\phi(t)$ are the unknown currents and potentials.
\begin{figure*}[t!]
  \centering
  \subfloat[]{\includegraphics[width=0.35\textwidth]{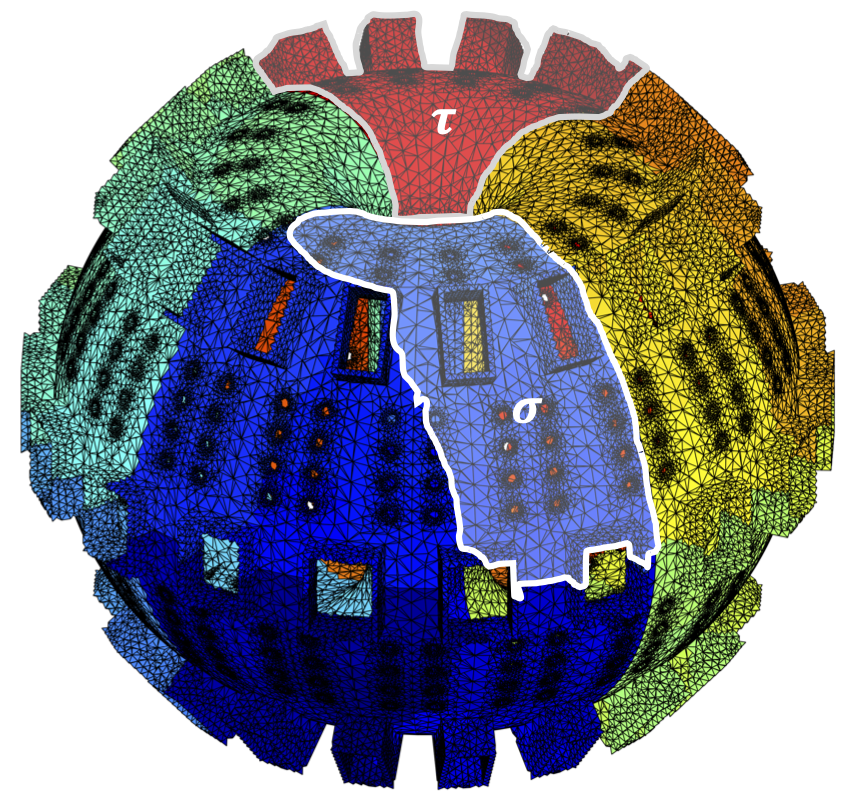}}
  \hfil
  \subfloat[]{\includegraphics[width=0.58\textwidth]{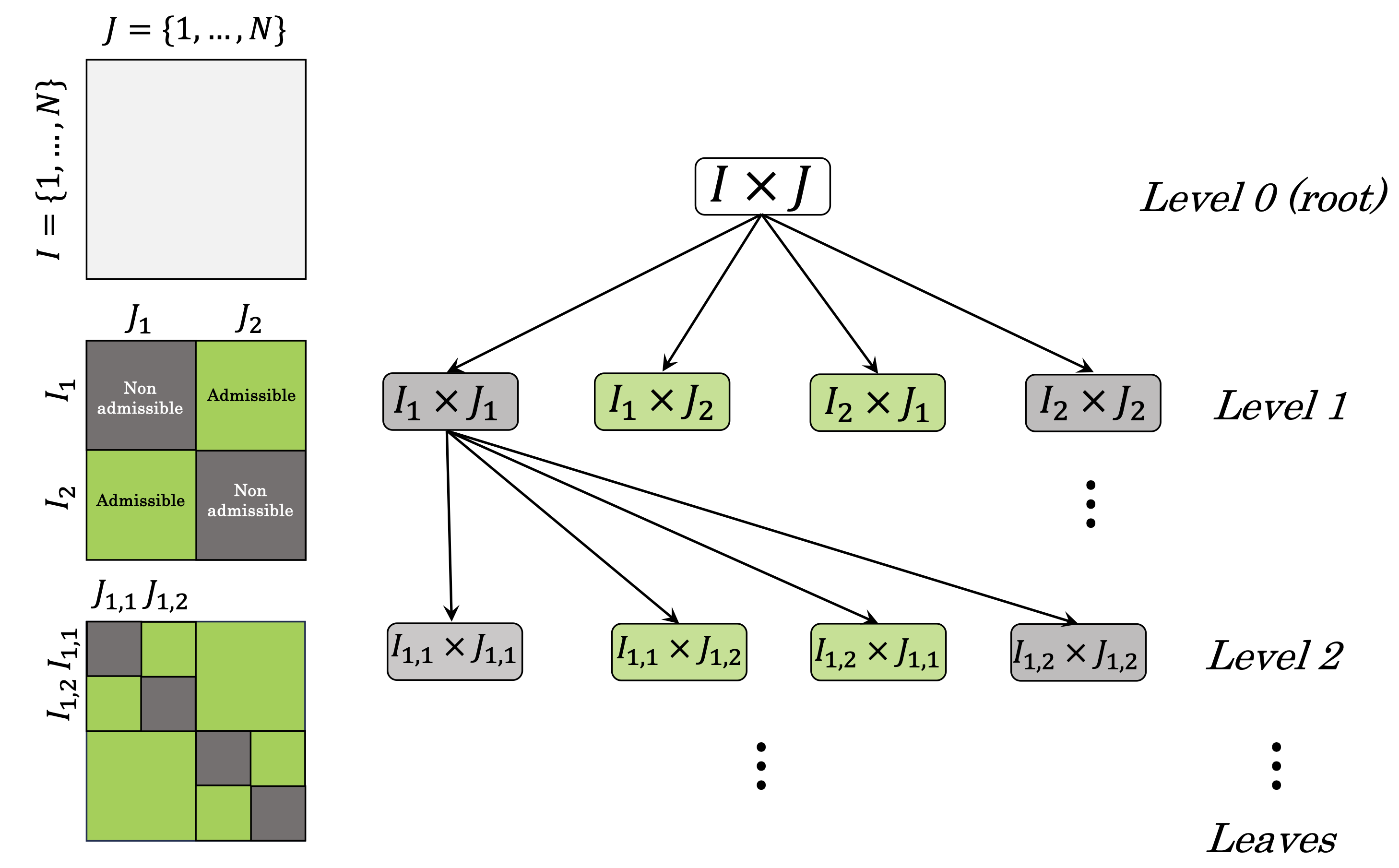}}
  \caption{Steps for construction of $\mathcal{H}$-matrix. (a) Reordering of DoFs in preprocessing. (b) Construction of block cluster tree of the product set $I\times J$ for the identification of low-rank (admissible) blocks.}
  \label{fig.construction_Hmat}
\end{figure*}
By applying Galerking testing to \eqref{eq.efie_div}, the following Differential Algebraic Equation (DAE) system is assembled:
\begin{equation}\label{eq.DAE}
\begin{bmatrix}
    \mathbf{L} & \mathbf{0} \\
    \mathbf{0} & \mathbf{0}
\end{bmatrix}
\frac{\text{d}}{\text{d}t}
\begin{bmatrix}
    \mathbf{j}(t) \\
    \boldsymbol{\Phi}(t)
\end{bmatrix}
=
-
\begin{bmatrix}
    \mathbf{R} & \mathbf{D}^\top \\
    \mathbf{D} & \mathbf{0}
\end{bmatrix}
\begin{bmatrix}
    \mathbf{j}(t) \\
    \boldsymbol{\Phi}(t)
\end{bmatrix}
-
\frac{\text{d}}{\text{d}t}
\begin{bmatrix}
    \mathbf{a_{ext}}(t) \\
    \mathbf{0}
\end{bmatrix},
\end{equation}
where the Degrees of Freedom (DoFs), namely the currents crossing the internal faces of the mesh and the potentials on each element, are collected in the arrays $\mathbf{j}$ and $\boldsymbol{\Phi}$. Matrices $\mathbf{L}$ and $\mathbf{R}$ are the inductance and resistance matrix of dimension $N_f\times N_f$, while matrix $\mathbf{D}$, of size  $N_v\times N_f$, is the incidence matrix of the equivalent circuit acting as the discrete counterpart of the divergence operator. Expressions of the entries of the aforementioned matrices can be found, e.g., in \cite{bettini2017volume}. 
In the context of magnetic confinement fusion devices, the magnet system composed of many current-carrying coils is responsible for the generation of the magnetic field required for the plasma current ramp-up, equilibrium control, and confinement of the charged particles~\cite{mitchell2008iter}. The external magnetic vector potential can be written as a linear map of the coil current $I(t)$. Thus, the last term in the right-hand side of \eqref{eq.DAE} can be written as:
\begin{equation}
    \frac{d\mathbf{a}_{ext}(t)}{dt}=\mathbf{B_i}\frac{dI(t)}{dt},
\end{equation}
where $\mathbf{B_i}$ identifies the current map.
It is worth pointing out that matrices $\mathbf{R}$ and $\mathbf{D}$ have a sparse pattern, while the inductance matrix $\mathbf{L}$, whose expression of the generic element $L_{ij}$ can be written as:
\begin{equation}
    L_{ij}=\mu_0\int_\Omega\int_\Omega G(\mathbf{r}_i,\mathbf{r}_j)\mathbf{w}_i\cdot\mathbf{w}_jd\Omega d\Omega,
\end{equation}
is a fully-populated matrix, due to the presence of the Green's kernel. The requirement for the memory storage of such matrix scales as $\mathcal{O}(N_f^2)$, rendering the cost prohibitive for realistic problems characterized by a significant number of unknowns. However, in the subsequent section, a solution to circumvent this challenge, leveraging the inherent smoothness property of the Green's function, is presented.


\subsection{Hierarchical matrix representation}\label{sub_sec.Hmatrix}

From the historical point of view, the main limitation of Integral Equation Methods (IEMs) for the solution of large-scale eddy current problems can be traced back to the fully populated pattern of the matrices arising from Galerkin's approach, as described in the previous section. Indeed, the storage complexity scales with $\mathcal{O}(N^2)$ while that for the inversion goes as $\mathcal{O}(N^3)$, where $N$ is the number of DoFs. Thanks to efficient techniques such as FFT-based approaches and data-sparse representations, the computational burden for the storage and manipulation of dense matrices arising from IEMs can be drastically reduced~\cite{alotto2015sparsification,bettini2020fast,cau2022fast,vacalebre2023low,10032268}. 
Among various compression techniques, hierarchical matrices, denoted as $\mathcal{H}$-matrices~\cite{voltolina2021optimized}, have emerged as a widely embraced approach, briefly elucidated in this section.

The $\mathcal{H}$-matrix approach relies upon the mathematical proprieties of the Green's kernel $G(\mathbf{r},\mathbf{r}^\prime)$, which exhibits a smooth behavior when $\mathbf{r}$ and $\mathbf{r}^\prime$ are far apart. This means that the interaction of distant source and target points is not only very weak but also similar. By generalizing this property, matrix blocks indexed by $\sigma$ and $\tau$, i.e., $\mathbf{L}(\sigma,\tau)$ resulting from the interaction of far elements exhibit a low-rank property and can be efficiently compressed. Here with the term compression, we mean that the matrix $\mathbf{L}_{\sigma,\tau}=\mathbf{L}(\sigma,\tau)$, of size $\#\sigma \times \#\tau$, where the symbol $\#$ refers to the cardinality of an index set, can be factorized, with a prescribed tolerance $\delta$, as:
\begin{equation}\label{eq.factorization}
    \mathbf{L}_{\sigma,\tau}\approx\mathbf{U_L}\mathbf{V_L}^\top,
\end{equation}
where $\mathbf{U_L}$ and $\mathbf{V_L}$ have dimension $\#\sigma\times r$ and $\#\tau\times r$, and $r$ is the effective rank of the factorization such that $r\ll\min\{\#\sigma,\#\tau\}$. This means that the storage cost of the matrix scales as $\mathcal{O}\left(r(\#\sigma+\#\tau)\right)$. 

With these concepts in mind, the primary objective in constructing the $\mathcal{H}$-matrix is to identify low-rank blocks suitable for compression. 
Firstly, a reordering of index sets is necessary to ensure that spatially distant points $\mathbf{r}_i$ and $\mathbf{r}_j$ correspond to matrix indices $i$ and $j$ that are similarly distant, as depicted in Fig.~\ref{fig.construction_Hmat}(a). The identification of low-rank blocks relies on leveraging the block cluster tree partition of the index sets associated with rows $I$ and columns $J$ of the matrix. Once the DoFs are appropriately reordered, the block cluster tree only necessitates the partitioning of index sets, commencing from the root $I\times J$, where $I,J=\{1,\ldots,N\}$. Each index set is then divided into two subsets. 
A matrix block becomes eligible for low-rank approximation if it satisfies the admissibility criterion~\cite{voltolina2019high}. Once this criterion is met, the blocks can be factorized in the form \eqref{eq.factorization}, leveraging the Adaptive Cross Approximation (ACA) technique~\cite{zhao2005adaptive}. This involves constructing matrices $\mathbf{U}_L$ and $\mathbf{V}_L$ through an adaptive procedure until a specified tolerance $\epsilon$ is achieved. Non-admissible blocks are partitioned iteratively until they reach a minimum user-defined block size. It's noteworthy that non-admissible blocks, stored in full form, are typically small, thereby ensuring a very high compression efficiency~\cite{voltolina2019high}.

\section{Linear elasticity model}\label{sec.FEM}
Despite the electromagnetic problem having a temporal dynamic, the mechanical problem is hereafter described by assuming static conditions. This is due to the different timescales driving the problem. Moreover, only small displacements are considered here, so a simple linear elasticity analysis can be set up. The governing Partial Differential Equation (PDE) of the structural problem is:
\begin{equation}\label{eq.linear_elasticity}
    -\boldsymbol{\nabla}\cdot \boldsymbol{\sigma}(\mathbf{U})=\mathbf{F},
\end{equation}
where $\boldsymbol{\sigma}$ is the Cauchy stress tensor and \textbf{F} the body force for unit volume. The Cauchy tensor is related to the strain-tensor $\boldsymbol{\epsilon}=0.5\left[\boldsymbol{\nabla}\mathbf{U}+\left( \boldsymbol{\nabla}\mathbf{U}\right)^\top\right]$ and the four-dimensional stiffness tensor $\mathcal{C}$ as $\boldsymbol{\sigma}=\mathcal{C}:\boldsymbol{\epsilon}$, where \textbf{U} is the displacement. 
Discretizing the computational domain with $N_v$ elements and exploiting a FEM approach, the PDE governing the linear elasticity problem is rewritten in the form:
\begin{equation}\label{eq.sys_structural}
    \mathbf{S}\mathbf{u}_m=\mathbf{f},
\end{equation}
where $\mathbf{S}$ is the stiffness matrix of dimension $3N_n\times 3N_n$, which is a function of the material properties, namely the density, Young's modulus, and Poisson's ratio. $N_n$ is the number of mesh nodes. The right-hand side vector $\mathbf{f}$ is the array of nodal external loads and $\mathbf{u}_m$ are the unknown nodal displacements. The expression for the entries of the stiffness matrix constructed with the 8-node brick elements can be found, e.g., in~\cite{Johnson2019brick}. The solution of the linear elasticity model requires the specification of the Dirichlet condition ($\mathbf{u}_{m,\Gamma_D}$=0) in portions of the boundary of the computational domain identified with $\Gamma_D$. The support structure of the VV of thermonuclear fusion devices is accurately designed to withstand the loads caused for example by electromagnetic, thermal, and seismic events during the operational life of the machine~\cite{ioki2009iter,martelli2021design}. The VV is usually supported with hinges distributed along the toroidal direction in the bottom part of each sector of the VV~\cite{martinez2014structural}. The bottom part of the support being fixed represents a Dirichlet boundary condition for the structural problem.
It is worth mentioning that, both the VV and its support structure are generally not axisymmetric, thus, the linear elasticity problem is always formulated as three-dimensional.

\section{Multiphysics coupling}\label{sec.coupling_VIE_FEM}
The unknown displacements of the static linear elastic problem are determined from the electromagnetic forces generated by the eddy currents and the external fields. 
The body force density at a given time is expressed as:
\begin{equation}  
\mathbf{F}(t,\mathbf{r})=\mathbf{J}(t,\mathbf{r})\times\mathbf{B}(t,\mathbf{r}),
\end{equation}
where $\mathbf{B}(t,\mathbf{r})$ is the magnetic flux density, sum of external and scattered fields due to the eddy currents, that is $\mathbf{B}(t,\mathbf{r})=\mathbf{B}_{eddy}(t,\mathbf{r})+\mathbf{B}_{ext}(t,\mathbf{r})$, and $\mathbf{J}(t,\mathbf{r})$ is the current density reconstructed from the solution of \eqref{eq.DAE} at a given time instant. It is worth noting that the body force density can be expressed with a combination of maps, indeed, the expansion \eqref{eq.j_phi_basis} defines a sparse map between the currents $\mathbf{j}$ and the components of the current density within each volume of the mesh, as:
\begin{align}
 \mathbf{J}_x(t,\mathbf{r})=&\mathbf{W}_x(\mathbf{r})\mathbf{j}(t), \nonumber \\  
 \mathbf{J}_y(t,\mathbf{r})=&\mathbf{W}_y(\mathbf{r})\mathbf{j}(t),\nonumber \\  
 \mathbf{J}_z(t,\mathbf{r})=&\mathbf{W}_z(\mathbf{r})\mathbf{j}(t) \label{eq.J_map_xyz}.
\end{align}
Matrices $\mathbf{W}_x(\mathbf{r}),\mathbf{W}_y(\mathbf{r}),\mathbf{W}_z(\mathbf{r})$ having a dimension of $N_v\times N_f$. The magnetic flux density vector due to the eddy current, instead, can be reconstructed from the current density field, assumed constant inside the mesh element $V$, through the analytical formula~\cite{fabbri2007magnetic}:
\begin{equation}\label{eq.fabbri_map}
    \mathbf{B}_{eddy}(t,\mathbf{r})=\frac{\mu_0}{4\pi}\sum_{S_f\in \partial V} \mathbf{J}(t)\times \mathbf{n}_fW_f(\mathbf{r}),
\end{equation}
where the sum is performed over all faces bounding the volume and the other symbols have the same meaning as in~\cite{fabbri2007magnetic}. By doing this, three dense matrices $\mathbf{K}_x(\mathbf{r}),\mathbf{K}_y(\mathbf{r}),\mathbf{K}_z(\mathbf{r})$, of dimension $N_v\times N_f$, are defined, allowing the construction of the components of $\mathbf{B}$ for each volume:
\begin{align}
 \mathbf{B}_{eddy,x}(t,\mathbf{r})=&\mathbf{K}_x(\mathbf{r})\mathbf{j}(t), \nonumber \\
 \mathbf{B}_{eddy,y}(t,\mathbf{r})=&\mathbf{K}_y(\mathbf{r})\mathbf{j}(t), \nonumber \\
 \mathbf{B}_{eddy,z}(t,\mathbf{r})=&\mathbf{K}_z(\mathbf{r})\mathbf{j}(t) \label{eq.B_map_xyz}.
\end{align}
It is worth noting that, the body force density can be fully constructed from the current DoFs $\mathbf{j}(t)$ through the defined matrices. Moreover, the dense matrices $\mathbf{K}_x,\mathbf{K}_y,\mathbf{K}_z$, exhibit low-rank properties, thus, they can be efficiently compressed with the previously described $\mathcal{H}$-matrices technique. Moreover, if the mesh is unchanged during the solution procedure, as is the case in our analyses, the aforementioned maps can be computed only once at the beginning of the procedure, thus, the cost for the reconstruction of the body force density is reduced to the cost of performing matrix-vector multiplications.

Using the maps \eqref{eq.J_map_xyz} and \eqref{eq.B_map_xyz} the right-hand side vector \textbf{f} for the linear elasticity model can be written as a function of the array of currents, solution of the EM problem: 
\begin{align}\label{eq.maps_j_to_f}
\mathbf{f}=\mathbf{P}\mathbf{F}=\mathbf{P}
\begin{bmatrix}
    \mathbf{F_x}  \\
    \mathbf{F_y} \\ \mathbf{F_z}
\end{bmatrix}
& =\mathbf{P}
\begin{bmatrix}
    \mathbf{J}_y\odot\mathbf{B}_z-\mathbf{J}_z\odot\mathbf{B}_y  \\
    \mathbf{J}_z\odot\mathbf{B}_x-\mathbf{J}_x\odot\mathbf{B}_z\\ \mathbf{J}_x\odot\mathbf{B}_y-\mathbf{J}_y\odot\mathbf{B}_x
\end{bmatrix} \nonumber \\
& =\mathbf{P}
\begin{bmatrix}
\mathbf{W}_y\mathbf{j}\odot(\mathbf{K}_z\mathbf{j}+\mathbf{B}_{ext,z})-\mathbf{W}_z\mathbf{j}\odot(\mathbf{K}_y\mathbf{j}+ \mathbf{B}_{ext,y})  \\
\mathbf{W}_z\mathbf{j}\odot(\mathbf{K}_x\mathbf{j}+\mathbf{B}_{ext,x})-\mathbf{W}_x\mathbf{j}\odot(\mathbf{K}_z\mathbf{j}+ \mathbf{B}_{ext,z}) \\
\mathbf{W}_x\mathbf{j}\odot(\mathbf{K}_y\mathbf{j}+\mathbf{B}_{ext,y})-\mathbf{W}_y\mathbf{j}\odot(\mathbf{K}_x\mathbf{j}+ \mathbf{B}_{ext,x}),
\end{bmatrix}
\end{align}
where the matrix \textbf{P} maps the force density \textbf{F} to the nodal loads \textbf{f}, and $\odot$ refers to the Hadamard product (i.e., element-wise product). The expression \eqref{eq.maps_j_to_f} reflects the tight connection between the EM and linear elasticity models, indeed, they share the same inputs which ultimately turn out to be the currents exciting the EM problem.

\section{Parametric model order reduction}\label{sec.pMOR}
In this section, the MOR of the electromagnetic and structural problems are described. 
As previously mentioned, the EM-ROM is first generated and then adopted also during the construction of the Structural-ROM. Different MOR strategies based on, e.g.,  Balanced Truncation, Moment Matching, or POD can be used to construct the ROMs \cite{benner2021model,10373844}.
Regardless of the adopted technique, MOR allows for projecting the original Full Order Model (FOM) (the EM or the structural one) onto a reduced order space. Thus, \eqref{eq.DAE} (the EM problem) and \eqref{eq.sys_structural} (the structural problem) are first rewritten in (descriptor) state space form:
\begin{align} \label{eq.SS_descriptor}
    \begin{split}
    \mathbf{{E}}\frac{d\mathbf{{x}}}{dt} &= \mathbf{{A}} \mathbf{{x}} +\mathbf{B_u} \mathbf{u}\\
            \mathbf{y} &= \mathbf{{C}} \mathbf{{x}} 
    \end{split}. 
\end{align}
It is worth noting that for the structural problem $\mathbf{E}=\mathbf{0}$. Then, the FOM \eqref{eq.SS_descriptor} is projected onto a reduced order state space model by using a projection matrix $\mathbf{V}$, i.e.,
\begin{align} \label{eq.ROM}
    \begin{split}
    \mathbf{\hat{E}}\frac{d\mathbf{\hat{x}}}{dt} &= \mathbf{\hat{A}} \mathbf{\hat{x}} +\mathbf{\hat{B}_u} \mathbf{u}\\
            \mathbf{\tilde{y}} &= \mathbf{\hat{C}} \mathbf{\hat{x}} 
    \end{split}, 
\end{align}
where $\mathbf{\hat{E}}=\mathbf{V^\ast\mathbf{E}\mathbf{V}}$ and similarly for  $\mathbf{\hat{A}}$, $\mathbf{\hat{B}_u}$, and $\mathbf{\hat{C}}$, while 
$\mathbf{\tilde{y}}\approx\mathbf{y}$. The superscript $\ast$ refers to the complex conjugate operator.
$\mathbf{V}$ is the projection matrix constructed by the adopted MOR strategy. The projection matrix $\mathbf{V}$ has dimension $N\times N_r$, where $N_r$ is the dimension of the reduced order space and $N_r \ll N$. 
Defying with $\eta$ the user-prescribed accuracy of the ROM, the projection matrix \textbf{V} is updated until the stop criterion is reached when $\epsilon<\eta$, with:
\begin{equation}\label{eq.error_ROM}
    \epsilon=\frac{\left\lVert [s\mathbf{E}-\mathbf{A}] \mathbf{\tilde{x}}-\mathbf{B_u u}\right\rVert}{\lVert\mathbf{B_u u}\rVert},
\end{equation}
where $\mathbf{\tilde{x}}$ is the approximate solution of the FOM, obtained as $\mathbf{\tilde{x}}=\mathbf{V}\mathbf{\hat{x}}$.
Since $\mathbf{V}$ has in general a limited number of columns (i.e., a small reduced order space is sufficient to accurately represent the dynamic of the quantity of interest stored in $\mathbf{y}$), the computational cost of solving the ROM \eqref{eq.ROM}  is much smaller than the one required to solve the FOM, making it compatible with real-time execution. 
 
Section~\ref{sub_sec.ROM_EM} deals with the steps for generating the ROM of the EM problem, while section~\ref{sub_sec.ROM_structural} described the guidelines for Structural-ROM construction. 

\subsection{Electromagnetic Reduced Order Model}\label{sub_sec.ROM_EM} 

In the following, a total of $N_{coils}$ are considered.
For the construction of the ROM of the EM problem, the original DAE in time domain \eqref{eq.DAE} is transformed in: 
\begin{equation}\label{eq.DAELap}
s
\underbrace{
\begin{bmatrix}
    \mathbf{L} & \mathbf{0} \\
    \mathbf{0} & \mathbf{0}
\end{bmatrix}
}_{\mathbf{E}}
\underbrace{
\begin{bmatrix}
    \mathbf{j} \\
    {\boldsymbol{\Phi}}
\end{bmatrix}
}_{\mathbf{x}}
=
\underbrace{
-
\begin{bmatrix}
    \mathbf{R} & \mathbf{D}^\top \\
    \mathbf{D} & \mathbf{0}
\end{bmatrix}
}_{\mathbf{A}}
\begin{bmatrix}
    {\mathbf{j}} \\
    {\boldsymbol{\Phi}}
\end{bmatrix}
\underbrace{
-
\begin{bmatrix}
    \mathbf{B_i}\\
    \mathbf{0}
\end{bmatrix}
}_{\mathbf{B_u}}
\underbrace{
sI
}_{\mathbf{u}}
\end{equation}
where $s\in \mathbb{C}$, and $\mathbf{j}$, $\boldsymbol{\Phi}$ are the unknown ``states'', collected in the array $\mathbf{x}$=[$\mathbf{j}$,$\boldsymbol{\Phi}$]$^\top$. 

By exploiting the linearity of the EM problem, the construction of the EM-ROM for a single excitation coil 
is repeated for each coil, resulting in $N_{coils}$ EM-ROMs.
The matrix to project the EM-FOM onto the Reduced Order space, i.e., $\mathbf{V}_{em,i}$, where $i$ indicates the $i$th coil, is constructed by using a standard POD method~\cite{benner2021model}.  
By using $\mathbf{V}_{em,i}$, with $i=1,\cdots,N_{coils}$, \eqref{eq.DAELap} is projected onto the reduced order space as in \eqref{eq.ROM}. The EM-ROMs can be then used to efficiently perform time-domain simulations in a very short computation time w.r.t. the time required to solve the EM-FOM.

\subsection{Structural Reduced Order Model}\label{sub_sec.ROM_structural}
Once the EM-ROMs for each excitation coil are constructed, the Structural-ROM can be finally generated. 
Differently from the EM-ROMs, which are independently constructed for any excitation coil by exploiting the linearity of the EM problem, a unique ROM is built for the structural problem and the effect of many EM sources is embedded within the MOR procedure itself. This is required since the forcing term of the Structural problem (i.e., the electromagnetic forces) has a non-linear dependence w.r.t. EM quantities.
The state space model for the structural problem is therefore written as:
\begin{equation} \label{eq.mec}
\mathbf{0}=\mathbf{S}\mathbf{u}_m-\mathbf{f}, 
\end{equation} 
where the right-hand side vector \textbf{f} is written according to \eqref{eq.maps_j_to_f}.
From \eqref{eq.maps_j_to_f}, \textbf{f} can be seen as a function of the coils' current and their time derivatives. 

Indeed, the EM load on the structure is due to the total induced current density and magnetic flux density, which can be obtained as a linear combination of the $\mathbf{J}_i$ and $\mathbf{B}_i$ related to each excitation coil, with $i=1,\ldots,N_{coils}$.
It is worth remembering that, since the EM-ROMs have already been generated, the solution to the EM problem and the time for constructing the force density through the maps \eqref{eq.J_map_xyz} and \eqref{eq.B_map_xyz} is drastically reduced.

The matrix to project the Structural-FOM into the Reduced Order space, i.e., $\mathbf{V}_{m}$, is constructed by using a standard POD method where the parameter space defined by the coils' currents and  $s$ in \eqref{eq.DAELap} are spanned. 

To generate the right-hand side snapshots of the Structural-FOM, time-domain simulations of the EM-ROMs are performed by varying the coils' currents and their dynamics according to realistic scenarios, to construct a suitable search space for the POD running on the Structural problem. From these time domain simulations, the current density vector $\mathbf{J}=[\mathbf{J}_x,\mathbf{J}_y,\mathbf{J}_z]$ and the magnetic flux density vector $\mathbf{B}=[\mathbf{B}_x,\mathbf{B}_y,\mathbf{B}_z]$  can be efficiently computed for each time step and for each mesh element by exploiting the linearity of the EM problem which has been projected onto a reduced order space, i.e.,
\begin{equation} \label{eq.Jtimestep}
\mathbf{J}^{k}_x=\mathbf{W}_x\sum^{N_{coils}}_{i}\mathbf{V}_{em,i}\mathbf{\hat{x}}^k_i
\end{equation}
where $k$ indicates the $k$th time step and $\mathbf{\hat{x}}^k_i$ is the solution of the EM-ROM for the $i$th coil at the $k$th time step. It is worth noting that an abuse of notation has been adopted in  \eqref{eq.Jtimestep} since $\mathbf{V}_{em,i}$ should be restricted only to the rows related to $\mathbf{j}$.
Similar equations can be written for $\mathbf{J}_y$, $\mathbf{J}_z$, $\mathbf{B}_x$, $\mathbf{B}_y$, and $\mathbf{B}_z$. 
It is worth noting that in \eqref{eq.Jtimestep}, the projection matrix $\mathbf{W}_x$ can be moved inside the sum. Thus, since $\mathbf{V}_{em,i}$ has dimension $N\times N_r$, with $N_r \ll N$, it is convenient  to compute $\mathbf{W}_x\mathbf{V}_{em,i}$ in advance to reduce the computation time. Once $\mathbf{J}^k_x$, $\mathbf{J}^k_y$, $\mathbf{J}^k_z$, $\mathbf{B}^k_x$, $\mathbf{B}^k_y$, and $\mathbf{B}^k_z$ have been constructed for each time step $k$, they are used by the POD algorithm to construct the projection matrix $\mathbf{V}_m$,  used to project \eqref{eq.mec} onto the reduced order space. 

\section{Solution of the time-dependent problem}\label{sec.time_integral}
Once the EM and structural ROMs are generated, the time-dependent analysis in the interval $[0,T]$ can be resolved with a limited computational effort, compatible with real-time (or faster than real-time) execution, which is essential for the monitoring of critical quantities. The time discretization is resolved by exploiting, for example, the $\theta$ method. Within each time step, the linear DAE \eqref{eq.ROM}, is solved for each input (coil) as:
\begin{equation}\label{eq.theta_method}
    \left[\theta \mathbf{\hat{A}}+\frac{1}{\tau}\mathbf{\hat{E}}\right]
    \hat{\mathbf{x}}_k=-\left[(1-\theta)\mathbf{\hat{A}}-\frac{1}{\tau}\mathbf{\hat{E}}\right]\hat{\mathbf{x}}_{k-1} 
    +\theta\mathbf{\hat{B}_u}u_k+(1-\theta)\mathbf{\hat{B}_u}u_{k-1},
\end{equation}
where $\theta\in[0,1]$, while $\tau$ is the width of time step. Then, at each $k$th time step, $\mathbf{J}^k_x$, $\mathbf{J}^k_y$, $\mathbf{J}^k_z$, $\mathbf{B}^k_x$, $\mathbf{B}^k_y$, and $\mathbf{B}^k_z$ are computed as in \eqref{eq.Jtimestep}. Then, $\mathbf{f}^k$ is computed as in \eqref{eq.maps_j_to_f} and the rhs of the Structural-ROM is obtained as:
\begin{equation} \label{eq.Vf}
\hat{\mathbf{f}}^k=\mathbf{V}^*_m\mathbf{f}^k.
\end{equation}
The procedure is repeated until the final integration time $T$ is reached.

The computational bottleneck of the approach is the computation of \eqref{eq.Jtimestep}, \eqref{eq.maps_j_to_f}, and \eqref{eq.Vf}. Indeed, such operations require projecting the Reduced Order solution onto the Full Order space. However, as previously discussed, the computational cost of  \eqref{eq.Jtimestep} can be effectively reduced by computing in advance $\mathbf{W}_x\mathbf{V}_{em,i}$. The computational cost of \eqref{eq.maps_j_to_f} is not dramatic since it requires to perform Hadamard products. Finally,  \eqref{eq.Vf} is a matrix-vector product but one dimension of $\mathbf{V}_m$ is the size of Structural-ROM, thus even computing \eqref{eq.Vf} in real-time is not costly. 
Numerical results reported in the next section show that the overall computational cost is compatible with real-time (or even faster than real-time execution) in standard hardware. 

To further reduce the computation time related to \eqref{eq.Jtimestep}, \eqref{eq.maps_j_to_f}, and \eqref{eq.Vf} a Discrete Empirical Interpolation Method (DEIM) approach can be used to avoid projecting the reduced order solution onto the full order space. This approach has not been implemented in this work since, for the considered problems, the computational cost was already compatible with the real-time execution of the model. 
However, \ref{Sec.Appendix}, illustrates how to apply DEIM to further speed up the computation time by trading a small amount of accuracy.

\section{Numerical example}\label{sec.numeric_ex}
The code is written in the MATLAB\textsuperscript{\textregistered} language and assisted with parallel MEX–FORTRAN functions based on OpenMP libraries for the fast evaluation of the matrix entries. The hierarchical matrix, in the Hierarchical Adaptive Low Rank (HALR) format~\cite{massei2022hierarchical}, has been constructed with the hm-toolbox, a MATLAB\textsuperscript{\textregistered} library accessible online through the following link:~\url{https://github.com/numpi/hm-toolbox}. 
The procedure outlined in the preceding sections is utilized for the real-time evaluation of structural displacement experienced by the VV of a thermonuclear fusion device during a VDE of the plasma current within the time interval of $[0,0.63]$ seconds.

\begin{figure*}[t!]
  \centering
  \subfloat[]{\includegraphics[width=0.47\textwidth]{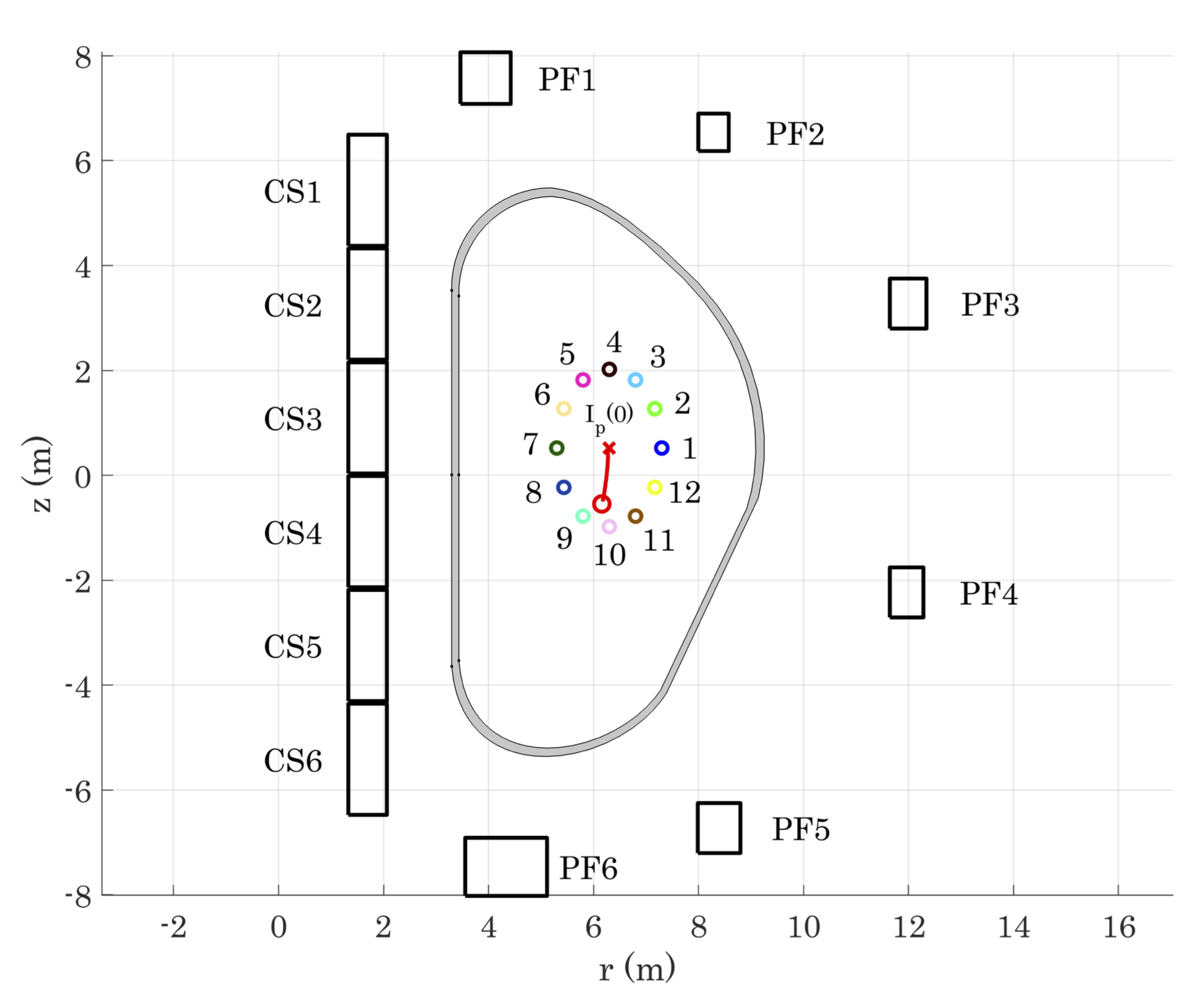}}
  \hfil
  \subfloat[]{\includegraphics[width=0.51\textwidth]{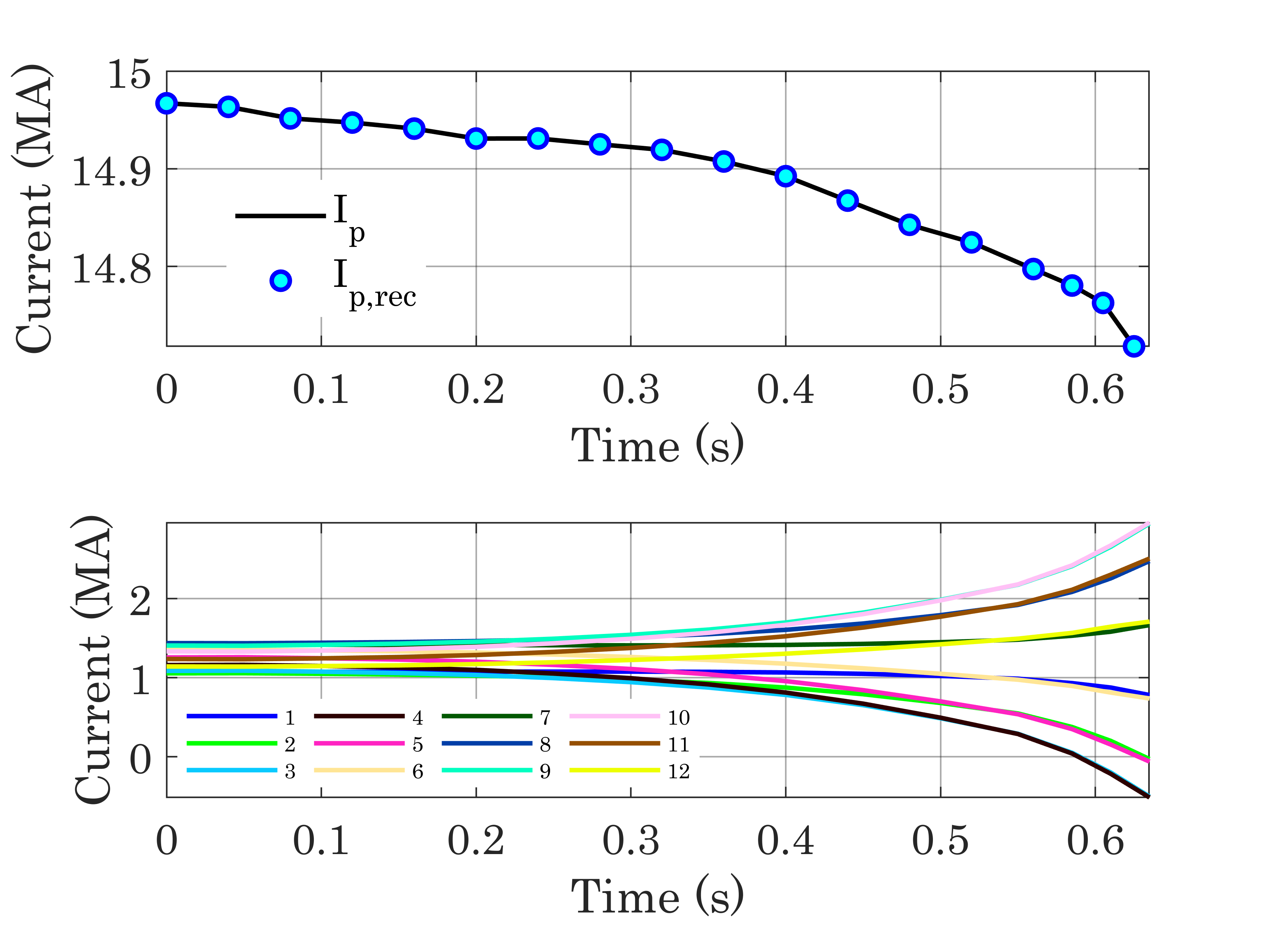}}
\caption{(a) Cross section of axisymmetric VV. Set of CS coils, and PF coils. The position of the plasma centroid at $t=0$ is marked with the red cross, while that at $t=0.63$ s with the red circle. The trajectory is also drawn. The equivalent current loops mimicking the plasma centroid movement are labeled with circular marks from 1 to 12. (b) The trend of plasma current $I_p(t)$ and that of equivalent current loops.} \label{fig.VDE_CAD}
\end{figure*}

The VDE results in a time-dependent plasma current $I_p(t)$ whose centroid position, under axisymmetric hypothesis, $\mathbf{r}_p(t)=\left(r_p(t),z_p(t)\right)$ also changes with time.
This behavior is common to most of the VDE observed in Tokamak plasmas. To mimic the spatial movement of the centroid $\mathbf{r}_p(t)$ with time, a set of $N_{eq}$ equivalent current loops, with a fixed position, are constructed around the initial centroid location $\mathbf{r}_p(t=0)$. The trend of each current flowing in the loops $I_{eq,i}(t)$ with $i=1,\ldots,N_{eq}$, is adjusted in such a way that the value of $I_p(t)$ at the position $\mathbf{r}_p(t)$ is correctly reconstructed until the position of the centroid lays within the circular crown defined by the $N_{eq}$ equivalent current loops. 
The described procedure is illustrated in Fig.~\ref{fig.VDE_CAD}(a), where the set of equivalent current loops and the position of the centroid at the initial time are highlighted, while in Fig.~\ref{fig.VDE_CAD}(b) the trend of the equivalent loop currents $I_{eq,i}(t)$ and the comparison between $I_p(t)$ and its reconstruction is given.

\begin{figure*}[t!]
  \centering
  \subfloat[]{\includegraphics[width=0.4\textwidth]{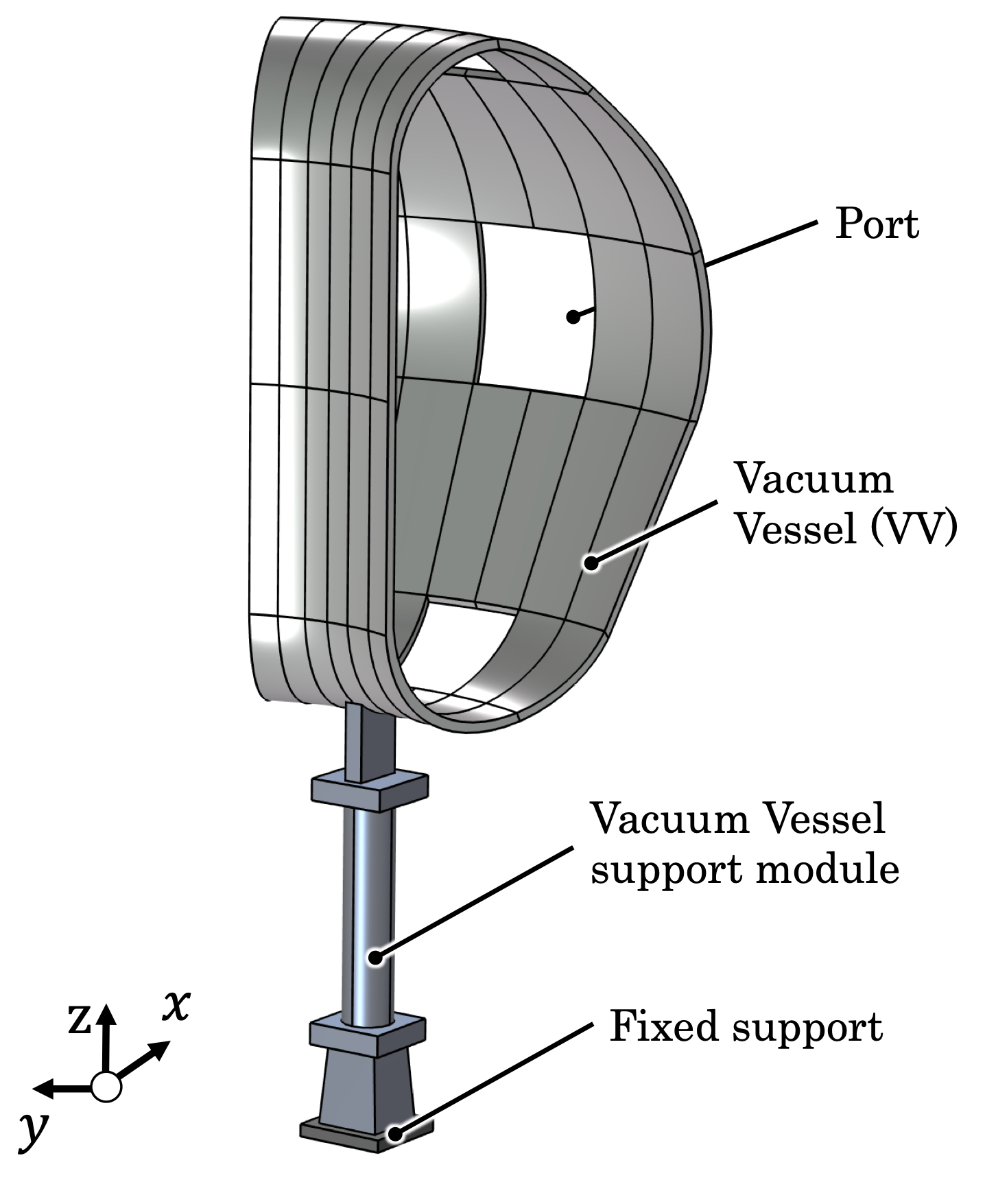}}
  \hfil
  \subfloat[]{\includegraphics[width=0.4\textwidth]{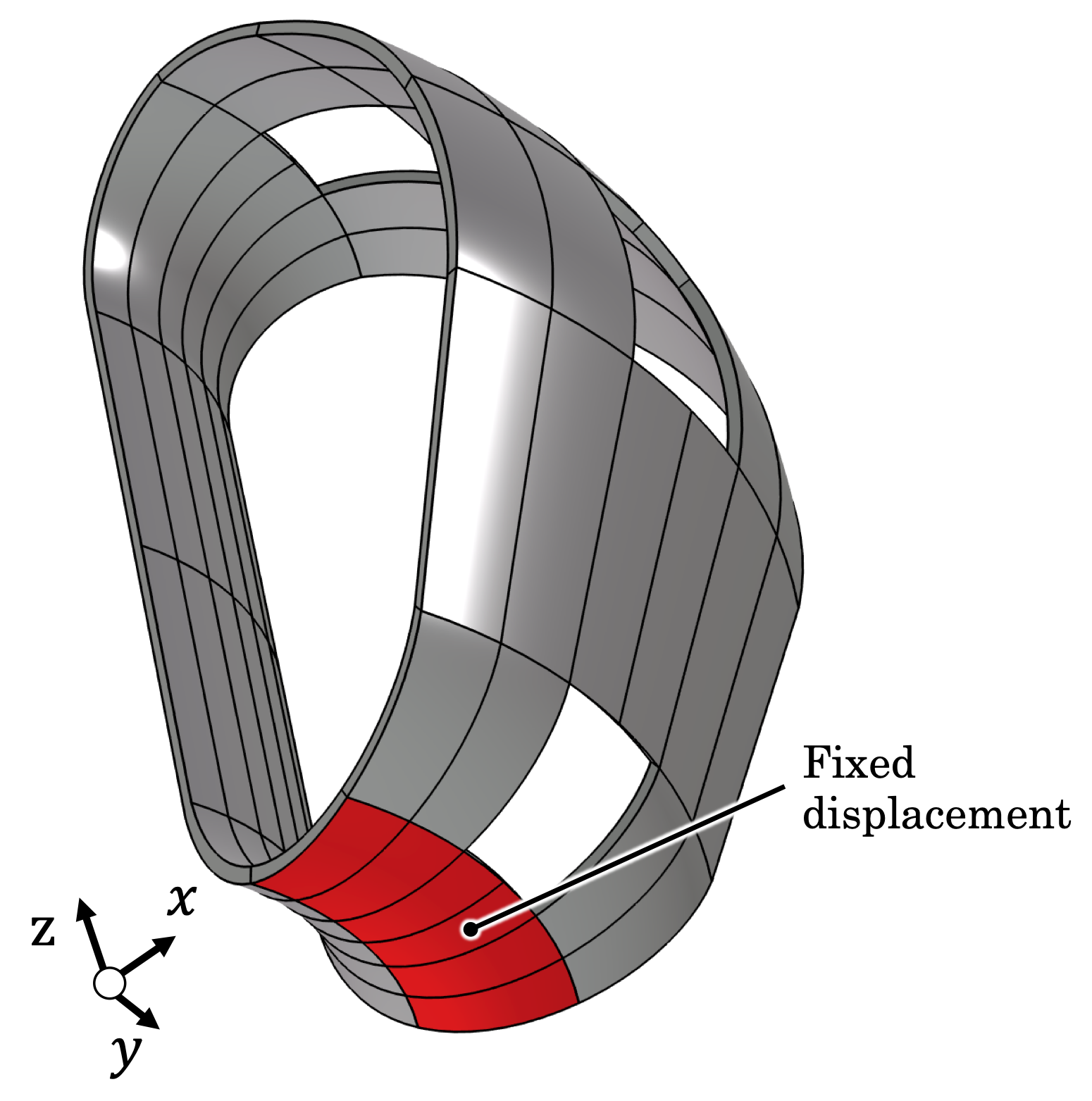}}
\caption{(a) Example of D-shaped VV sector support module. (b) Fixed displacement area.} 
\label{fig.VV_support}
\end{figure*}

Following this approach, a set of current-carrying coils with prescribed time-dependent current, on which the ROM can be construed, are defined. Thus, following the notation used in section~\ref{sec.pMOR},  $N_{eq}=N_{coils}$. It is worth mentioning that the construction of such equivalent currents is not strictly necessary, instead, a sequence of coils following the path of the centroid can be generally used. However, the latter approach results in a more complicated analysis. For simplicity, due to the VDE, the currents in the Central Solenoid (CS) and Poloidal Field (PF) coils are kept constant with the values listed in Table~\ref{tab.CS_PFcurr}. This assumption is justified by the different time scales involving the phenomena~\cite{xu2016analytical}. Thus, the CS and PF systems are only used to set up a background static magnetic flux density $\mathbf{B}_0$ which does not affect the eddy current problem but enters in the evaluation of force density. 
\begin{table}[t!]
  \centering
  \caption{CS and PF coil currents.}
  \label{tab.CS_PFcurr}
  \begin{tabular}{ccc}
    \hline
    Coil numbering & $I_{CS}$ (MA) & $I_{PF}$ (MA) \\
    \hline
    1 & -4.8 & 5.7 \\
    2 & -4 & -2.8 \\
    3 & -20.8 & -5.8 \\ 
    4 & -20.8 &-4.8  \\ 
    5 & -10 & -7.8 \\ 
    6 & 5.2 & 16.8  
\end{tabular}
\end{table}
In the following numerical examples, the described approach for VDE modelling is applied to understand the structural response of the D-shaped VV of thermonuclear fusion machines. Note that in general, the VV structure is subjected to loads of different nature, however, in what follows, only the Lorentz force is considered, so that, the effect of thermal loads and gravity are neglected. Concerning the linear elasticity model, despite the VV support structure for which an example is illustrated in Fig.~\ref{fig.VV_support} (a), are intrinsically 3D, in the following a simplified assumption is made. Without losing generality, the Dirichlet conditions are imposed in the points on the surface highlighted in red in Fig.~\ref{fig.VV_support} (b). By doing this, an axisymmetric fixed support directly attached to the VV is defined. 

The material parameters used for the VV, required by the EM and structural analyses are summarized in Table~\ref{tab.VV_material}.
\begin{table}[t!]
  \centering
  \caption{Material properties of VV.}
  \label{tab.VV_material}
  \begin{tabular}{cc}
    \hline
    Property & Value \\
    \hline
    Electric resistivity & 7.4$\cdot10^{-7}$ ($\Omega$m) \\
    Young's modulus & 193 (GPa) \\
    Poisson's ratio & 0.25 \\ 
    Density & 8000 (kg/m$^3$)   
\end{tabular}
\end{table}
The example given in section~\ref{sub_sec:VV_axi} validates the approach proposed in this paper against an axisymmetric model resolved within the commercial FEM software COMSOL{\textsuperscript{\textregistered}} Multiphysics. However, for a realistic VV of a thermonuclear fusion device, due to many ports providing access for heating, diagnostics, and remote handling operations, the axisymmetric assumption is too rough, and 3D models should be used instead. 
For this reason, the second numerical example described in section~\ref{sub_sec:VV_3D} shows the capabilities of the proposed approach in determining the displacement of a 3D VV similar to that of ITER Tokamak.
 
\subsection{D-shaped Vacuum Vessel}\label{sub_sec:VV_axi}
The first numerical example consists of an axisymmetric conductor resembling a D-shaped VV of a thermonuclear fusion device. The geometry is constructed by revolving the cross-section of the VV reported in Fig.~\ref{fig.VV_support} (a) neglecting the ports. The VV is characterized by a major radius of 6.2~m, a height of approximately 11~m, and a thickness of 14~cm. 
The effectiveness of the resulting EM-Structural ROM is tested against the commercial FEM software COMSOL{\textsuperscript{\textregistered}} Multiphysics. The EM field is excited with the CS and PF currents listed in Table~\ref{tab.CS_PFcurr}, while the VDE is modeled with the approach described above by defining the set of $N_{eq}=12$ equivalent plasma currents whose trend is given in Fig.~\ref{fig.VDE_CAD} (b). 
The time integration lasts for $T=0.63$~s. Within COMSOL{\textsuperscript{\textregistered}}, the EM-structural analysis is performed by using the time-domain \emph{Magnetic Fields} (mf) module and, at each time step, the stationary \emph{Solid Mechanics} (solid) module. The two modules are coupled via the Lorentz force computation. 

\begin{figure}[t!]
    \centering
    \includegraphics[width=0.6\linewidth]{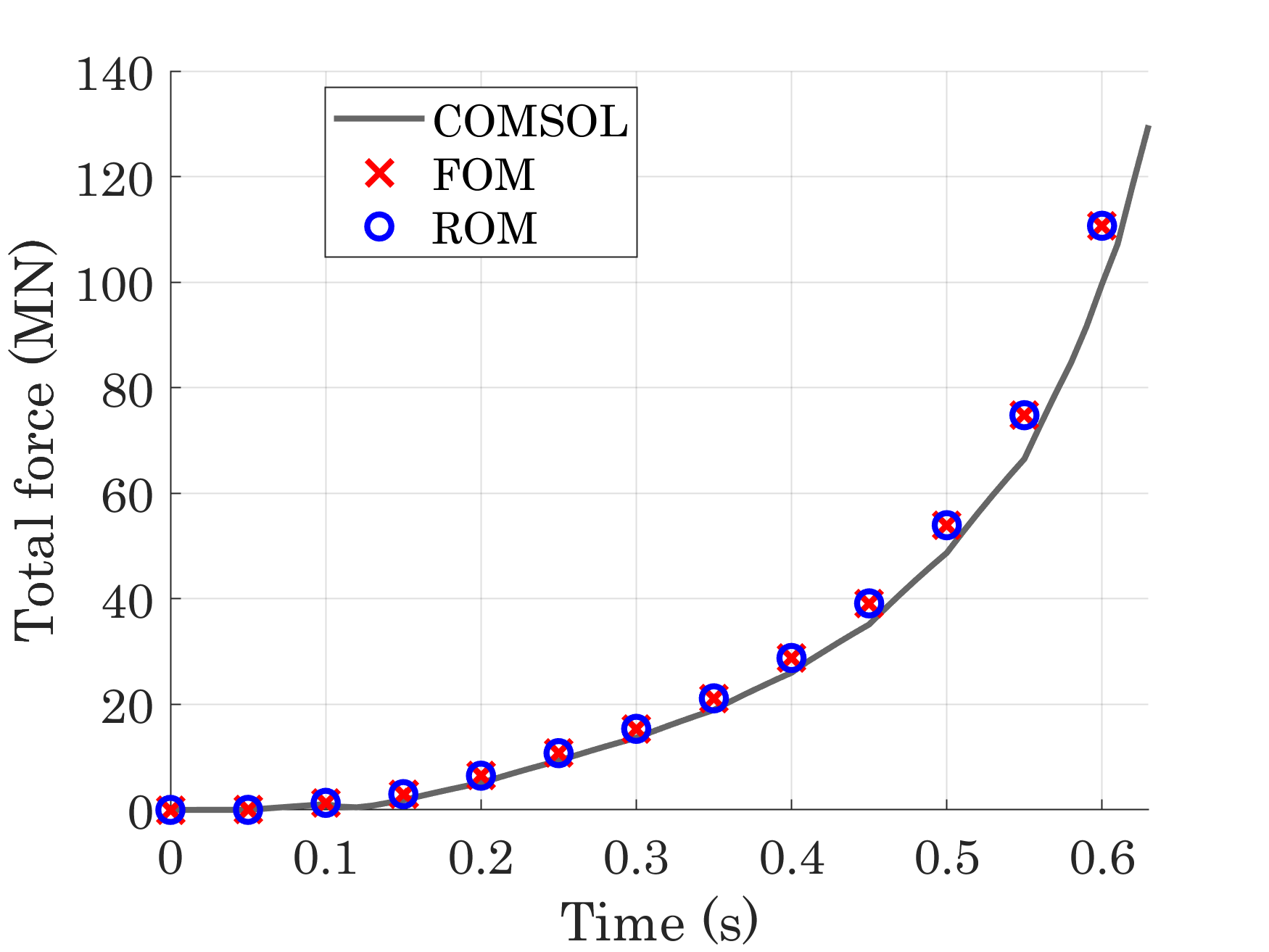}
    \caption{Trend of total force over the VV during the VDE evaluated with the 2D axisymmetric model in COMSOL{\textsuperscript{\textregistered}}, the FOM, and the ROM.}
    \label{fig.2Daxial_force}
\end{figure}

For the EM-structural ROM, the VV is discretized with 22k tetrahedral elements, resulting in $N_f=60k$ current DoFs. The storage of the full inductance matrix \textbf{L} would require $\approx 29$~GB, however, thanks to the data-sparse representation explained in section~\ref{sub_sec.Hmatrix}, $\approx 7$~GB are needed, resulting in a compression ratio of 76\%. The EM-ROMs are constructed with a tolerance equal to $\eta=10^{-3}$, 
resulting in $N_{eq}=12$ EM-ROMs which are stored with low memory usage being the greatest projection matrix $\mathbf{V}_{em,max}$ of dimension $N \times 5$. This means that the reduced matrices $\mathbf{\hat{E}}$ and $\mathbf{\hat{A}}$ of each EM-ROM, are not greater than $5 \times 5$. 
For the Structural-ROM, a tolerance  $\eta=10^{-3}$ is also set,
and  a projection matrix $\mathbf{V}_{m}$ of dimension $3N_n \times 27$ is obtained. 
As depicted in Fig.~\ref{fig.2Daxial_force}, the integration of the force density \textbf{F} across the entire volume of the VV obtained from the ROM overlaps the one obtained from the corresponding FOM. Moreover, results are in excellent agreement with the ones obtained from the axisymmetric COMSOL{\textsuperscript{\textregistered}} model.

\begin{figure*}[t!]
  \centering
  \subfloat[]{\includegraphics[width=0.35\textwidth]{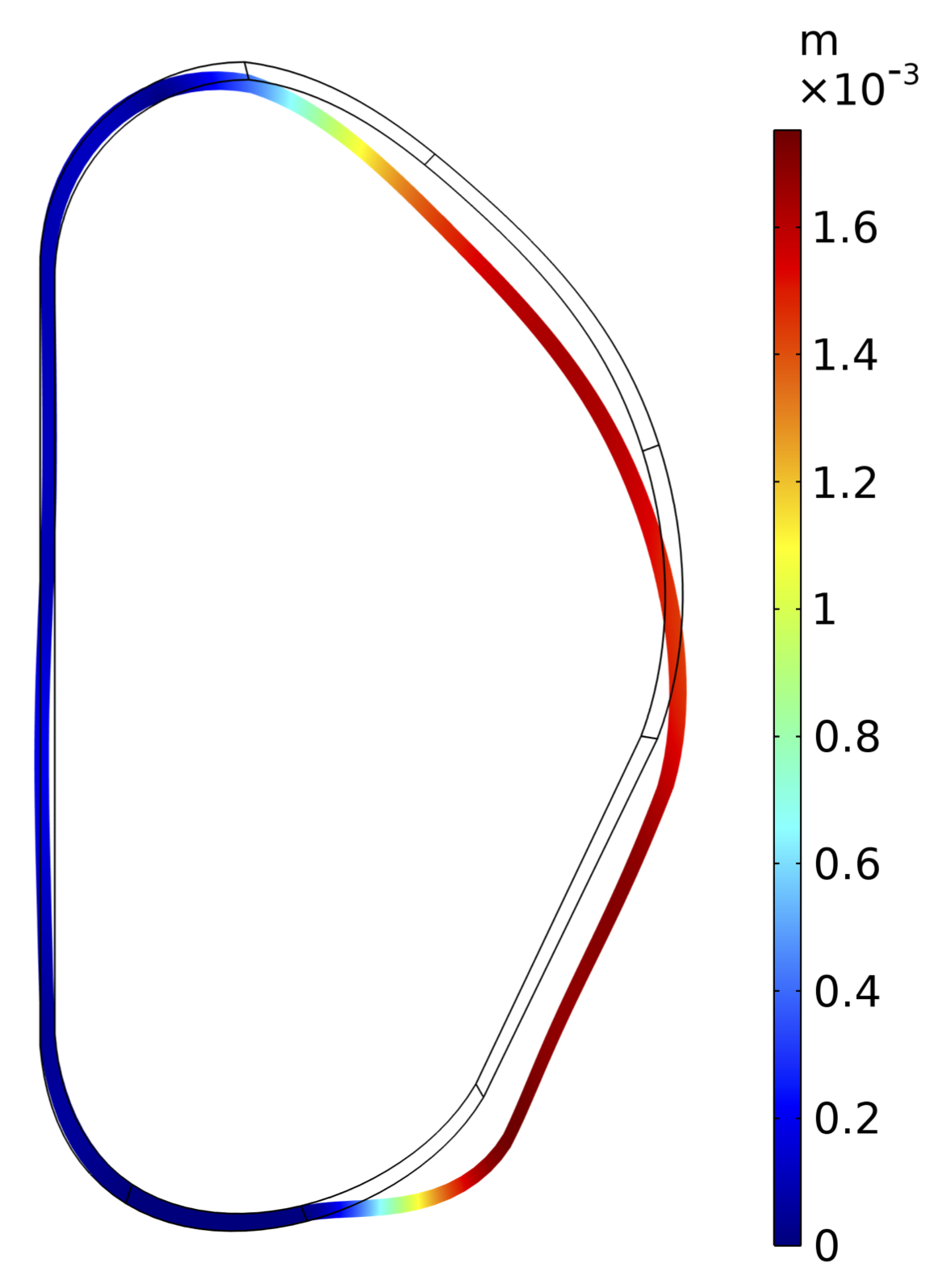}}
  \hfil
  \subfloat[]{\includegraphics[width=0.36\textwidth]{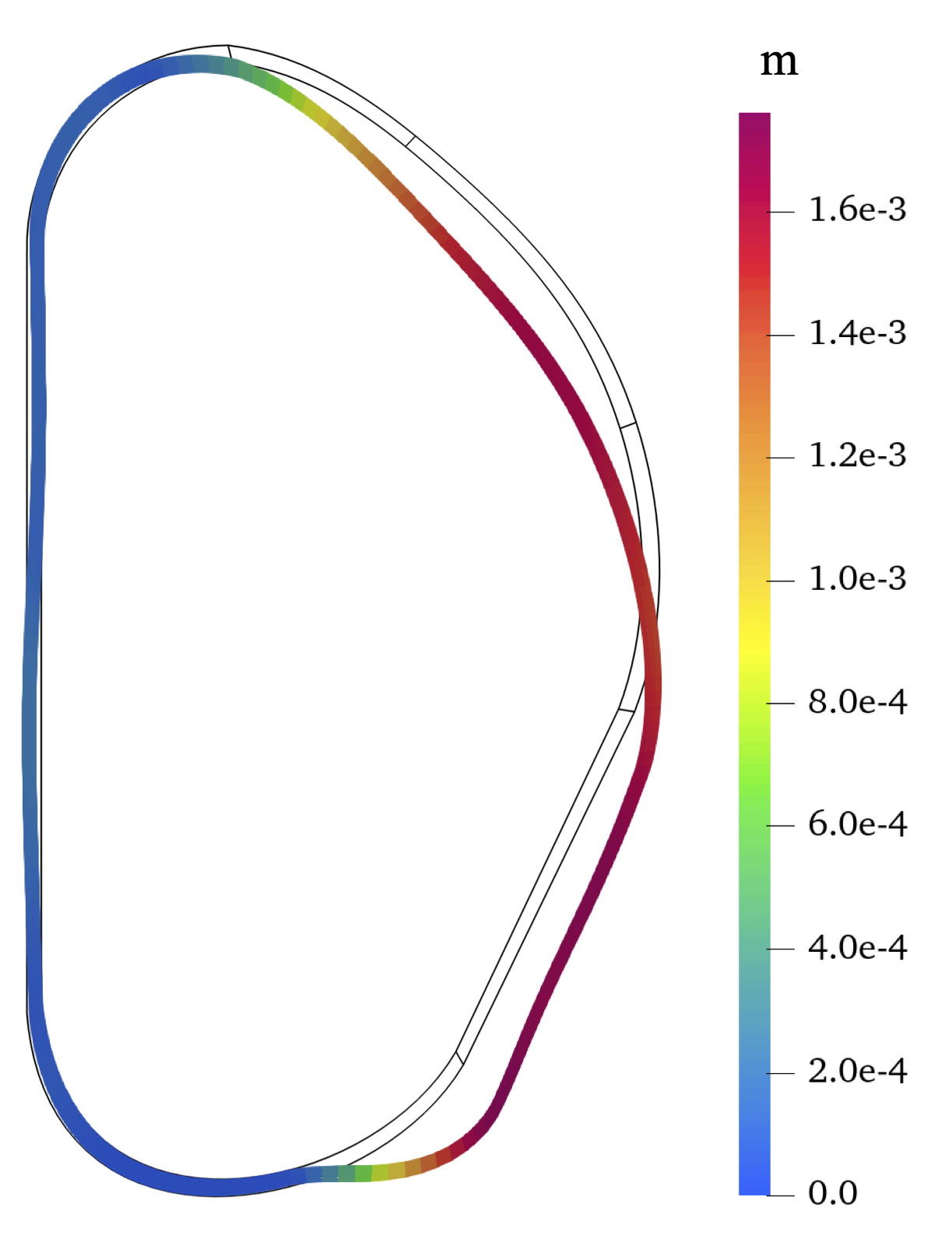}}
\caption{(a) Norm of the displacement field at $t=0.63$ s evaluated with the 2D axisymmetric model in COMSOL{\textsuperscript{\textregistered}}. (b) The magnitude of the displacement field at $t=0.63$ s evaluated with the ROM. Graphically the displacement is amplified by a factor $\alpha=300$.} 
  \label{fig.2Daxial_displace}
\end{figure*}

\subsection{ITER-like Vacuum Vessel}\label{sub_sec:VV_3D}
As an additional example, let's examine the influence of a VDE on a VV design reminiscent of that found in the ITER Tokamak.
The whole VV is constructed by fully revolving the VV sector of 30° illustrated in Fig.~\ref{fig.VV_support}(a). It is worth noting that, due to the presence of the ports, the problem is not axisymmetric, thus, a full 3D analysis is mandatory. The VDE is modeled with the approach described in section~\ref{sec.numeric_ex}. The EM and structural ROMs are constructed in approximately 4~h and 45~min respectively.  
Nonetheless, it is noteworthy that the maximum number of bases constituting the projection matrices $\mathbf{V}_{em}$ and $\mathbf{V}_m$ is confined to be 5 and 25, respectively, obtained with the same tolerance values of the previous example. 

\begin{figure}[t!]
    \centering
    \includegraphics[width=0.8\linewidth]{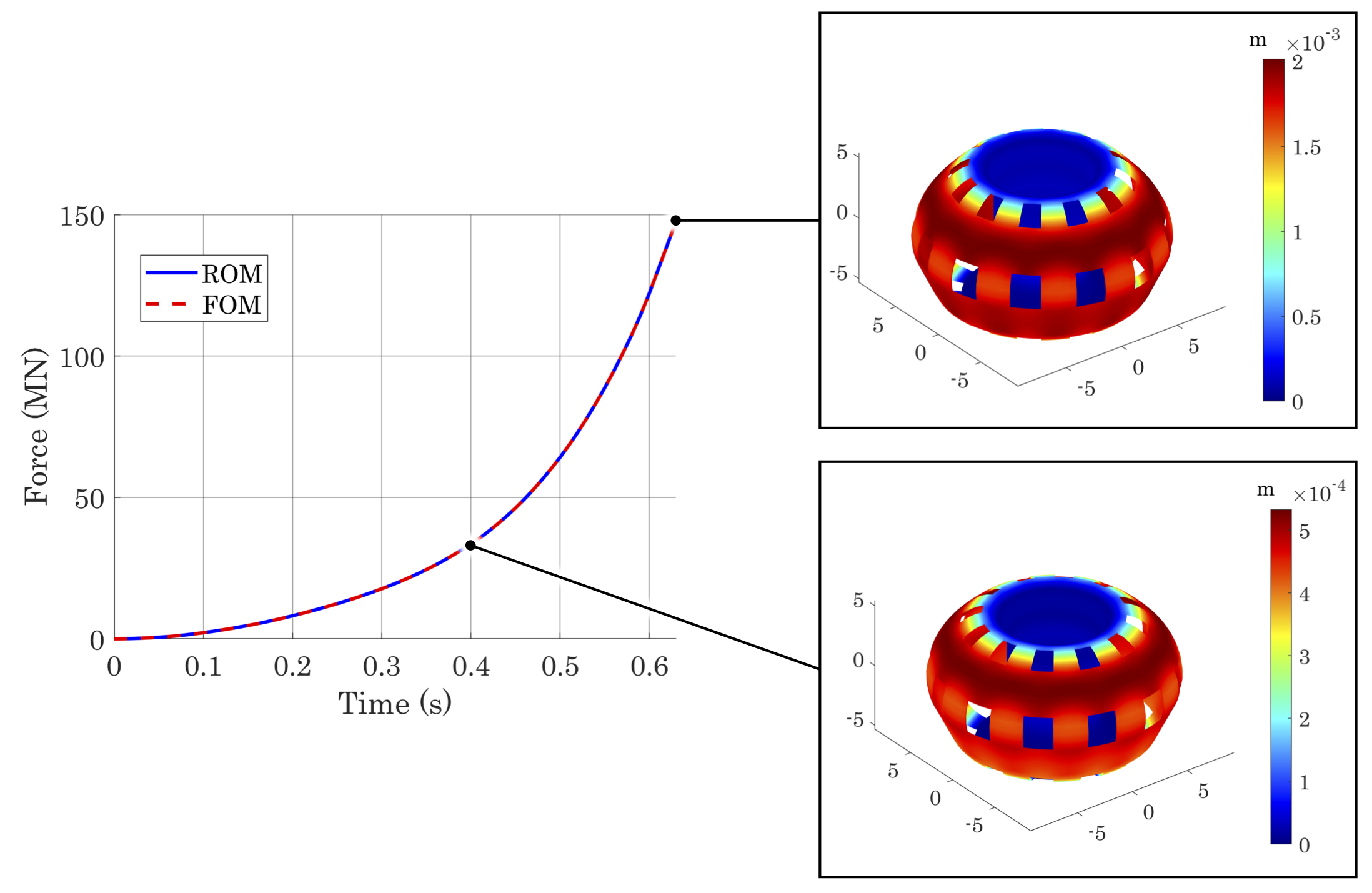}
    \caption{Trend of total force over the ITER-like Vacuum Vessel during the VDE evaluated with the FOM, and the ROM. The displacement magnitude is also reported at $t=0.4$ s and at $t=0.63$ s. Graphically the displacement is amplified by a factor $\alpha=600$.}
    \label{fig.3D_force}
\end{figure}

Again, to test the validity of the methodology described in this paper, the total force in the structure over the time window is evaluated both with the ROM and the FOM. The trend of the force illustrated in Fig.~\ref{fig.3D_force} shows a good agreement between the two approaches. The magnitude of the displacement is also reported in Fig.~\ref{fig.3D_force} at two time instants, corresponding to $t=0.4$ s and at $t=0.63$ s.
The time needed to solve the multiphysics EM-structural problem varies significantly between the ROM and the FOM. Specifically, the time integration with the ROM through \eqref{eq.theta_method} demands less than 0.5~s, demonstrating ``faster than real-time" claims, whereas the FOM's time integration takes approximately 90~h. This discrepancy underscores the impracticality of employing the FOM for real-time monitoring of structural displacement in large-scale conductive structures while demonstrating the efficiency of the proposed approach.

\section{Conclusions}\label{sec.conclusion}
The analysis of the structural deformation of conductive structures in magnetic confinement fusion devices due to VDE is of great interest for the operation of current experimental and future machines. In this paper, a MOR approach for the real-time estimation of the structural deformations has been presented. The MOR
is essential for reducing the computational burden imposed by the FOM in establishing a Digital Twin of the machine. The MOR approach allows the construction of EM and structural ROMs which can be used to solve the transient EM-Structural analysis in real-time (or faster than real-time). The proposed methodology is designed to accurately assess the structural deformation caused by electromagnetic forces acting upon a realistic 3D VV of thermonuclear fusion devices during a VDE. 

It is worth mentioning that the EM problem for the construction of the EM-ROM, as discussed in this paper, has been addressed using the VIE method, accelerated by hierarchical matrices ($\mathcal{H}$-matrices). However, the same procedure can also be implemented using a FEM approach.

The numerical results illustrate the outstanding performance of the MOR procedure in accelerating the transient analysis of VDE, making it compatible with real-time simulation of the VDE phenomenon and allowing ``faster than real-time" highly accurate predictions suitable for incorporation in control algorithms.

\appendix
\section{Discrete Empirical Interpolation Method} \label{Sec.Appendix}
For the sake of simplicity, but without loosing generality, in this section it is assumed that $\mathbf{J}_y$ and $\mathbf{B}_z$ are nonzero, while the other components of the current density vector and the magnetic flux density are zero. The right-hand side of the structural problem in the reduced order space is obtained as
\begin{equation} \label{eq.fVPF}
\hat{\mathbf{f}}=\mathbf{V}_m^*\mathbf{P}\mathbf{F}=\mathbf{V}_m^*\mathbf{P}(\mathbf{J}_y\odot\mathbf{B}_z).
\end{equation}
For the sake of simplicity, a slight abuse of notation is applied in \eqref{eq.fVPF} since $\mathbf{P}$ and $\mathbf{V}_m$ should be restricted to the columns related to the $x$ component of the force. During the real-time computation of \eqref{eq.fVPF}, $\mathbf{J}_y$ and  $\mathbf{B}_z$ (from now on, the sub-scripts are omitted to simplify the notation) are obtained as in \eqref{eq.Jtimestep}, i.e., 
$\mathbf{J}=\mathbf{W}\mathbf{V}_{em}\hat{\mathbf{x}}$ and $\mathbf{B}=\mathbf{K}\mathbf{V}_{em}\hat{\mathbf{x}}$, where $\hat{\mathbf{x}}$ is the reduced order solution of the EM problem at a given time step and only one coil is assumed again for the sake of simplicity. 
Thus, \eqref{eq.fVPF} can be rewritten as:
\begin{equation} \label{eq.fVPW}
\hat{\mathbf{f}}=\mathbf{V}_m^*\mathbf{P}\big((\mathbf{W}\mathbf{V}_{em}\hat{\mathbf{x}})\odot(\mathbf{K}\mathbf{V}_{em}\hat{\mathbf{x}})\big). 
\end{equation}
The problem of evaluating \eqref{eq.fVPW} in real-time is that the solution of the reduced order EM problem is first projected onto the full order space where the force is computed. Then, the force is projected onto the reduced order space of the structural problem. Although only Hadamard products are performed in the full order space, when the size of the FOMs is very large, \eqref{eq.fVPW}  may become incompatible with real-time execution. To solve this problem, the DEIM approach can be used. The interested reader is referred to~\cite{chaturantabut2010nonlinear} for more information about DEIM. At the same time, hereinafter it is shown how DEIM applied to \eqref{eq.fVPW} allows avoiding to work on the full order space.
DEIM allows for approximating \eqref{eq.fVPW} as:
\begin{equation}
\hat{\mathbf{f}} \approx   \mathbf{V}^*_m\mathbf{P}\mathbf{Z}  \Big[(\mathbf{S}\mathbf{Z})\backslash\big((\mathbf{S}\mathbf{W}\mathbf{V}_{em}\mathbf{x})\odot(\mathbf{S}\mathbf{K}\mathbf{V}_{em}\mathbf{x})\big)\Big],
\end{equation}
where:
\begin{itemize}
    \item $\mathbf{Z}$ is a $N \times k$ matrix with the DEIM bases. It can be obtained, for instance, by performing a Truncated Singular Value Decomposition (TSVD) of a matrix containing snapshots of $(\mathbf{W}\mathbf{V}_{em}\hat{\mathbf{x}})\odot(\mathbf{K}\mathbf{V}_{em}\hat{\mathbf{x}})$. $N$ is the dimension of the full order model while $k$ is order chosen for the TSVD, with $k\ll N$,
    \item $\mathbf{S}$ is a  $k \times N$ selection matrix that projects the full-order solution onto a subset of interpolation points selected by the DEIM procedure, i.e., by using MATLAB notation, $\mathbf{S}=\mathbf{I}(\mathbf{s},:)$, where $\mathbf{I}$ is the identity matrix and $\mathbf{s}$ is the array of indexes of the interpolation points. 
\end{itemize}

To reduce the computational cost of \eqref{eq.fVPW}, the following considerations hold:
\begin{itemize}
    \item $\mathbf{V}^*_m\mathbf{P}\mathbf{Z}$ is pre-calculated leading to a small $N_r \times k$ matrix,
    \item $\mathbf{S}\mathbf{Z}$ is pre-calculated leading to a small $k \times k$ matrix,
    \item $\mathbf{S}\mathbf{Z}$ is pre-calculated leading to a small $k \times k$ matrix,
    \item $\mathbf{S}\mathbf{W}\mathbf{V}_{em}\mathbf{x}$ and $\mathbf{S}\mathbf{K}\mathbf{V}_{em}\mathbf{x}$ are pre-calculated leading to small $N_r \times k$ matrices. 
\end{itemize}
By doing that, the computational cost of \eqref{eq.fVPW} no longer depends on the size of the FOM and the computational cost is reduced to $\mathcal{O}(N_r k + k^2)$, making it compatible with real-time implementation even for large-scale problems. 

\bibliographystyle{elsarticle-num}
\bibliography{bibliography.bib}

\end{document}